\documentclass[11pt]{amsart}
\usepackage[english]{babel}
\usepackage{amsmath}
\usepackage{amssymb}
\usepackage{amsthm}
\usepackage{libertine}
\usepackage[all]{xy}
 \newtheorem{thm}{Theorem}[section]
 \newtheorem{cor}[thm]{Corollary}
 
 \newtheorem{prop}[thm]{Proposition}
  \newtheorem{pb}[thm]{Problem}
 \theoremstyle{definition}
 \newtheorem{defn}[thm]{Definition}
 \theoremstyle{remark}
 \newtheorem{rem}[thm]{Remark}
 
 \numberwithin{equation}{section}

\newcommand{\scal}[1]{\left<#1\right>}
\newcommand{\Hq}{\mathbb H}
\newcommand{\Sq}{\mathbb S}

\newcommand{\N}{\mathbb{N}}

\newcommand{\R}{\mathbb{R}}      

\newcommand{\C}{\mathbb{C}}

\title[Fock and Hardy spaces: Clifford Appell case]{Fock and Hardy spaces: Clifford Appell case }

\author[D. Alpay]{Daniel Alpay}
\address{(DA) Schmid College of Science and Technology, Chapman University, Orange 92866, CA, US}
\email{alpay@chapman.edu}
\thanks{Daniel Alpay thanks the Foster G. and Mary McGaw Professorship in
Mathematical Sciences, which supported this research.}
\author[K. Diki]{Kamal Diki}
\thanks{Kamal Diki thanks Chapman university for kind hospitality during the period in which a part of this paper was written. This research is supported by the project INdAM Doctoral Programme in Mathematics and/or Applications Cofunded by Marie Sklodowska-Curie Actions, acronym: INdAM-DP-COFUND-2015, grant number: 713485. }
\address{(KD) Politecnico di
Milano\\Dipartimento di Matematica\\Via E. Bonardi, 9\\20133 Milano,
Italy}
\email{kamal.diki@polimi.it}

\author[I. Sabadini]{Irene Sabadini}
\address{(IS) Politecnico di
Milano\\Dipartimento di Matematica\\Via E. Bonardi, 9\\20133 Milano\\Italy}
\email{irene.sabadini@polimi.it}

\begin{document}
\maketitle
\begin{abstract}
In this paper, we study a specific system of Clifford-Appell polynomials
and in particular their product. Moreover, we introduce a new family of
quaternionic reproducing kernel Hilbert spaces in the framework of
Fueter regular functions. The construction is based on a general idea
which allows to obtain various function spaces, by specifying a suitable
sequence of real numbers. We focus on the Fock and Hardy cases in this
setting, and we study the action of the Fueter mapping and its range.
 \end{abstract}


\section{Introduction}
A set of polynomials $\lbrace P_n \rbrace_{n\in\mathbb{N}}$ satisfying an identity with respect to the real derivative that takes $P_n$ to $nP_{n-1} $ is called an Appell system \cite{Appell1880}. In the classical case, where $x$ is interpreted as a real or complex variable,  the standard monomials $P_n(x)=x^n$ form an Appell set, but also the famous Hermite, Bernoulli and Euler polynomials are examples of Appell sets. The importance of such polynomials in various settings is well known, and we mention here, with no pretense of completeness their relevance in probability theory and stochastic process since they can be connected to random variables, see \cite{Bao2015}, they were used also to study optimal stopping problems related to Lévy process in \cite{Salm2011}.

  Moving to the hypercomplex analysis setting, namely analysis for functions with values in a Clifford algebra, in particular quaternions, we have various functions theories, associated with different differential operators.
  In this paper we will treat the quaternionic case.\\
  In the slice hyerholomorphic setting, Appell systems can be obtained by simply extending the variable in use to become  hypercomplex, and so we have that, for example, the standard monomials in the  quaternionic variable are among them with respect to the slice derivative.

  But these sets of polynomials were studied also in the setting of quaternionic and Clifford analysis with respect to the hypercomplex derivative, see  \cite{CMF2017,CMF2011,DKS2019, MF2007,Pena2011}. It turns out that the Appell systems in this framework play a similar role  as the complex monomials do to define elementary functions in terms of their power series like cosine, sine, exponential, etc. This fact opens a variety of questions also in relation to various function spaces including Fock, Hardy, Bergman, Dirichlet spaces, etc. Moreover, various questions arise about their associated operators such as creation, annihilation, shift and backward shift operators. Some different operators related to Fock spaces in the Clifford setting were considered also in \cite{CFR2011}. What makes Appell systems in quaternionic and Clifford analysis rather peculiar, is the fact that the function theory has been developed using the so-called Fueter polynomials, see  \cite{BDS}, \cite{GHS}, and these polynomials do not satisfy the Appell property in general. However, a series expansion for hyperholomorphic functions is possible using both the approaches. \\
   In order to define and study quaternionic reproducing kernel Hilbert spaces
   the approach that makes use of the Appell systems looks very promising and allows to define the associated operators. We will show that using a special set of Clifford Appell polynomials, denoted by $\{Q_n\}$, we can introduce various functions spaces denoted by $\mathcal{HM}_b$ whose elements are converging series of the form $\sum Q_n a_n$, where the quaternionic coefficients $a_n$ satisfy suitable conditions which depend on a given sequence $b=(b_n)$ of real (in fact rational) numbers.
    This approach is rather general, and it is used also in the slice hyperholomorphic setting in which the series under consideration are of the form $\sum q^n a_n$, where $q$ denotes the quaternionic variable and give rise to spaces denoted by $\mathcal{HS}_c$, $c=(c_n)$.
    \\
    In this paper we treat the case of the quaternionic Fock and the Hardy spaces which have been already studied in the slice setting but are new in the Fueter regular framework combined with the Appell polynomials. For this reason, these spaces are called Clifford-Appell Fock space and Clifford-Appell Hardy space, respectively.
   \\
   One problem of the system $\{Q_n\}$ is that if we multiply two such polynomials we do no obtain an element in the system. This is expected provided the non-commutative setting and in fact hyperholomorphic functions can be multiplied using the so-called CK-product. With the polynomials $Q_n$ there is the additional problem of remaining within the Appell system and in fact we show how this can be achieved. This technical result opens the possibility to prove several results and also to introduce creation, annihilation and shift operators.
    \\
    An advantage of our description is that we can prove that the function spaces $\mathcal{HM}_b$ and $\mathcal{HS}_c$ for suitable choices of $b,c$, can be related using the Fueter mapping theorem.
   \\
   The structure of the paper is the following: in Section 2 we revise notations and preliminary results that we need in the sequel. In Section 3 we introduce some quaternionic reproducing kernel Hilbert spaces (QRKHS) based on a specific Appell system, and prove different properties on such kind of polynomials. We show also that, under suitable conditions, any axially Fueter regular function can be expanded in terms of these Appell polynomials. In Section 4 we focus more on the Fock space in this setting. In particular, we study different properties related to the notions of creation, annihilation operators and Segal-Bargmann transforms. In Section 5 we treat the Hardy space case, and study different properties related to the shift and backward shift operators. Finally, in Section 6 we show how the Fueter mapping acts by sending spaces of slice hyperholomorphic functions into spaces of Fueter regular functions. Moreover, we show that in some special cases the Fueter mapping acts as an isometric isomorphism up to a constant.

\section{Preliminary results}
We recall some basic facts on quaternions and on the two sets of Cauchy-Fueter and slice hyperholomorphic functions.
The skew field of quaternions is defined to be
$$\Hq=\lbrace{q=x_0+x_1i+x_2j+x_3k\quad ; \ x_0,x_1,x_2,x_3\in\R}\rbrace$$ where the imaginary units satisfy the multiplication rules $$i^2=j^2=k^2=-1\quad \text{and}\quad ij=-ji=k, jk=-kj=i, ki=-ik=j.$$
The conjugate and the modulus of $q\in\Hq$ are defined by
$$\overline{q}=Re(q)-\vec{q} \quad \text{where} \quad Re(q)=x_0, \quad \vec{q}=x_1i+x_2j+x_3k$$
and $$\vert{q}\vert=\sqrt{q\overline{q}}=\sqrt{x_0^2+x_1^2+x_2^2+x_3^2},$$
 respectively.
Notice that the quaternionic conjugation satisfy the property $\overline{ pq }= \overline{q}\, \overline{p}$ for any $p,q\in \Hq$.
Moreover, the unit sphere $$\lbrace{\vec{q}=x_1i+x_2j+x_3k;\text{ } x_1^2+x_2^2+x_3^2=1}\rbrace$$ coincides with the set of all  imaginary units given by $$\mathbb{S}=\lbrace{q\in{\Hq};q^2=-1}\rbrace.$$
Sometimes we denote $e_1=i,$ $e_2=j$ and $e_3=k$.
\\
We recall the classical notion of Fueter regular functions also called "hyperholomorphic functions", for more details one can see \cite{CSSS2004, GHS}:
\begin{defn} Let $U\subset \Hq$ be an open set and $f:U\longrightarrow \Hq$ a real differentiable function. We say that $f$ is (left) Fueter regular or regular for short if $$\partial f(q):=\displaystyle\left(\frac{\partial}{\partial x_0}+i\frac{\partial}{\partial x_1}+j\frac{\partial}{\partial x_2}+k\frac{\partial}{\partial x_3}\right)f(q)=0, \forall q\in U.$$
\\
The quaternionic right linear space of Fueter regular functions is denoted by $\mathcal{R}(U)$.
\end{defn}

The right Fueter regular functions can be defined just by taking the imaginary units on the right of the derivatives of the function $f$. The quaternionic monomials $P_n(q)=q^n$ are not Fueter regular. However, there exist some other important functions in this theory, the so-called Fueter variables, defined by
\begin{equation}
\zeta_l(x)=x_l-e_lx_0, \textbf{  } l=1,2,3.
\end{equation}

These functions play the same role that complex monomials play in complex analysis. For example, a series expansion for Fueter regular functions is obtained using these Fueter variables. A suitable product that allows to preserve the regularity in this setting is the so-called C-K product, denoted $\odot$. Given two Fueter regular functions $f$ and $g$, we take their restriction to $x_0=0$ and consider their pointwise multiplication. Then, we take the Cauchy-Kowalevskaya extension of this pointwise product, which exists and is unique, to define $f\odot g$, see \cite{GHS}.

 A more recent theory of quaternionic regular functions was introduced  and studied in several directions during the last years, see for example \cite {ACS_book,CSS,CSS_book, GentiliSS},
namely the theory of slice hyperholomorphic functions that we recall briefly. In the definition below, for a fixed $I\in \Sq$, $\C_{I}=\R+I\R$ denotes  the complex plane whose variable is $q=x+Iy$, and we set $\Omega_I := \Omega \cap \C_I$.
\begin{defn}
A real differentiable function $f: \Omega \longrightarrow \Hq$, on a given domain $\Omega\subset \Hq$, is said to be a (left) slice hyperholomorphic function if, for very $I\in \Sq$, the restriction $f_I$ to $\C_{I}$, is holomorphic on $\Omega_I$, that is it has continuous partial derivatives with respect to $x$ and $y$ and the function
$\overline{\partial_I} f : \Omega_I \longrightarrow \Hq$ defined by
$$
\overline{\partial_I} f(x+Iy):=
\dfrac{1}{2}\left(\frac{\partial }{\partial x}+I\frac{\partial }{\partial y}\right)f_I(x+yI)
$$
vanishes identically on $\Omega_I$. The set of slice hyperholomorphic functions is denoted by $\mathcal{SR}(\Omega)$.
\end{defn}

The paper  \cite{AlpayColomboSabadini2014} studies the slice hyperholomorphic quaternionic Fock space $\mathcal{F}_{Slice}(\Hq)$,  defined for a given $I\in{\mathbb{S}}$  to be
$$\mathcal{F}_{Slice}(\Hq):=\left\lbrace{f\in{\mathcal{SR}(\Hq); \, \displaystyle  \frac{1}{\pi} \int_{\C_I}\vert{f_I(p)}\vert^2 e^{-\vert{p}\vert^2}d\lambda_I(p) <\infty}}\right\rbrace,$$
 where $f_I = f|_{\C_I}$ and $d\lambda_I(p)=dxdy$ for $p=x+yI$. The definition of this space does not depend on the choice of $I$. It was also proved that this quaternionic Fock space can be characterised in terms of the slice hyperholomorphic power series as follows $$\mathcal{F}_{Slice}(\Hq)=\left\lbrace{\displaystyle \sum_{k=0}^\infty q^ka_k; \textit{ } a_k\in \Hq \text{ :} \sum_{k=0}^\infty k!|a_k|^2<\infty}\right\rbrace.$$

 Its associated Segal-Bargmann transform was studied  in \cite{DG1.2017} by considering the slice hyperholomorphic kernel obtained making use of the normalized Hermite functions $(\eta_n)_{n\geq 0}$. The explicit expression of this kernel is given by
 \begin{equation}\label{BargKer}
  \mathcal{A}_\Hq^S(q,x):=\displaystyle \sum_{k=0}^\infty\frac{q^k}{\sqrt{k!}}\eta_k(x)=e^{-\frac{1}{2}(q^2+x^2)+\sqrt{2}qx}, \textbf{ } \forall (q,x)\in \Hq\times \R.
 \end{equation}

Then, for any quaternionic valued function $\varphi$ in $L^2(\R,\Hq)$ the slice hyperholomorphic Segal-Bargmann transform is defined by
\begin{equation}\label{BargTra}
\displaystyle \mathcal{B}_\Hq^S(\varphi)(q)=\int_\R \mathcal{A}_\Hq^S(q,x) \varphi(x)dx.
\end{equation}
In the same spirit different famous spaces of slice hyperholomorphic functions such as Hardy, Besov, Bloch, Dirichlet and Bergman spaces were studied in \cite{AlpayColomboSabadini2015, CastilloColomboGantnerCervantes2015, ColomboCervantesSabadini2015}.

\section{A new family of QRKHS of Fueter regular functions: General setting}
Let us consider the quaternionic polynomials defined by
\begin{equation} \label{Qk2}
\displaystyle Q_k(q)=\sum_{j=0}^k T^k_jq^{k-j}\overline{q}^j, \text{} q\in\Hq, \textbf{ } k\geq 0
\end{equation}

where
\begin{equation}
\displaystyle T^k_j:=\frac{k!}{(3)_k}\frac{(2)_{k-j}(1)_j}{(k-j)!j!}=\frac{2(k-j+1)}{(k+1)(k+2)}
\end{equation}
 and $(a)_n=a(a+1)...(a+n-1)$ is the Pochhammer symbol.
\begin{rem}
Notice that the polynomials $(Q_k)_{k\geq 0}$  given by \eqref{Qk2} are Fueter regular on $\Hq$. Moreover, they form an Appell system with respect to the hypercomplex derivative $\displaystyle\frac{\overline{\partial}}{2}$. i.e, for all $k\geq 1$ we have the Appell property
\begin{equation}
\displaystyle\frac{\overline{\partial}}{2}Q_k=kQ_{k-1}.
\end{equation}
 For more details on such properties of these polynomials one can consult for example \cite{CMF2011} and \cite{DKS2019}.
\end{rem}

 For $q\in \Hq,$ let
 \begin{equation}
  {\rm {Exp}}(q):=\sum_{k=0}^\infty\frac{Q_k(q)}{k!}
 \end{equation}
 be the generalized Fueter regular exponential function considered in the paper \cite{CMF2011}. Then, we introduce the following
\begin{defn}\label{DefHb}
Let $\Omega$ be a domain in $\Hq$. Let $c=(c_k)_{k\in\N}$  and $b=(b_k)_{k\in\N}$ be two non decreasing sequences with $c_0=b_0=1$. Then, associated to $b$ and $c$ we define
\begin{enumerate}
\item The subspace of Fueter regular functions defined by  $$\mathcal{HM}_{b}(\Omega)=\left\lbrace{\displaystyle \sum_{k=0}^\infty Q_k\alpha_k; \textit{ } \alpha_k\in \Hq \text{ : } \sum_{k=0}^\infty b_k|\alpha_k|^2<\infty}\right\rbrace.$$
\item The subspace of slice hyperholomorphic functions defined by $$\mathcal{HS}_{c}(\Omega)=\left\lbrace{\displaystyle \sum_{k=0}^\infty q^kf_k; \textit{ } f_k\in \Hq \text{ :  } \sum_{k=0}^\infty c_k|f_k|^2<\infty}\right\rbrace.$$

\end{enumerate}

\end{defn}
 Given $f=\displaystyle \sum_{k=0}^\infty Q_k\alpha_k \text{ and } g=\displaystyle \sum_{k=0}^\infty Q_k\beta_k$ in $\mathcal{HM}_b(\Omega)$ we define the Hermitian inner product given by $$\scal{f,g}_{\mathcal{H}_b}=\displaystyle\sum_{k=0}^\infty b_k\overline{\alpha_k}\beta_k.$$

\begin{rem}
We note that, by specifying the sequence $c$, $\mathcal{HS}_c$ include different spaces of slice hyperholomorphic functions such as Fock, Hardy, Dirichlet and generalized Fock spaces. Such spaces are the quaternionic counterpart of the complex version introduced in \cite{AlpayColomboSabadini2019}.
\end{rem}

 We are interested in two main problems in this setting:
\begin{pb}
 Study the counterparts of the spaces introduced in Definition \ref{DefHb} by suitably chosing the sequence $b$ in order to include in this framework of Cauchy-Fueter regularity : Fock, Bergman, Hardy, Dirichlet spaces, etc.
\end{pb}
 In this paper, we will treat the Fock and Hardy cases that correspond, respectively, to the sequences $b_k=k!$ and $b_k=1$, $\forall k \geq 0$.
\begin{pb}\label{pb2}
Study the range of the Fueter mapping on $\mathcal{HS}_c$ and see when it is possible to obtain spaces of regular functions of the form  $\mathcal{HM}_b$. More in general, we ask if using the Fueter mapping it is possible to get information on the sequence $(b_k)$ in terms of the given datum $(c_k)$ ?
\end{pb}

\begin{rem}
We note that the answer to Problem \ref{pb2} for Fock and Bergman cases were considered in \cite{DKS2019}. See also \cite{AlpayColomboSabadini2014, ColomboCervantesSabadini2015} for the slice hyperholomorphic setting.  The answer in these  two cases is given by:
\begin{enumerate}
\item The Fock case: $$c_k=k! \text{ and } b_k=\displaystyle \frac{k!}{(k+1)(k+2)}, \textbf{  } \forall k \geq 0.$$
\item The Bergman case: $$c_k=\displaystyle \frac{1}{k+1}\text{ and } b_k=\displaystyle \frac{1}{(k+1)^2(k+2)^2(k+3)}, \textbf{  } \forall k \geq 0.$$
\end{enumerate}
\end{rem}

 We will show that, under suitable conditions, for some special choices of the sequence $b$ in Definition \ref{DefHb} we have the estimate:

\begin{equation}
  \displaystyle |f(q)|\leq \left(\sum_{k=0}^\infty \frac{|q|^{2k}}{b_k}\right)^{\frac{1}{2}} \Vert f \Vert_{\mathcal{HM}_b}, \textbf{ } f\in\mathcal{HM}_b(\Omega), \textbf{ }q\in\Omega.
\end{equation}

In these cases, we can also prove that
 $\mathcal{HM}_b(\Omega)$ are right quaternionic reproducing kernel Hilbert spaces with reproducing kernel given by
 \begin{equation}
  K_{\mathcal{HM}_b(\Omega)}(q,p)=\displaystyle \sum_{k=0}^\infty \frac{Q_k(q)\overline{Q_k(p)}}{b_k}, \textbf{  } \forall (q,p)\in\Omega\times \Omega.
 \end{equation}

Furthermore, in such situations $\displaystyle\left\lbrace \frac{Q_k}{\sqrt{b_k}}\right\rbrace_{k\geq 0}$ form an orthonormal basis of $\mathcal{HM}_b(\Omega)$.

Now, we will prove an interesting result on the Appell polynomials $(Q_k)_{k\geq 0}$ useful to compute their C-K product.

\begin{prop}\label{QkCKQs}
 Let $k,s \geq 0$. Then, for any $q=x_0+\vec{q}\in\Hq$ we have

$$\displaystyle (Q_k\odot Q_s)(q)=\frac{c_kc_s}{c_{k+s}}Q_{k+s}(q), $$
where $\odot$ is the C-K product and $c_l:=\displaystyle \sum_{j=0}^l(-1)^jT_j^l, \textbf{  } \forall l\geq 0.$
\end{prop}
\begin{proof}
 Since $Q_k$ and $Q_s$ are Fueter regular functions on $\Hq$, their C-K product $Q_k\odot Q_s$ is also Fueter regular. Then, we use the formula of the C-K extension, see \cite[Theorem 11.38]{GHS}, given by $$CK[h(\vec{q}\,)](q)=\exp\left(-x_0\partial_{\vec{q}\,}\right)[h(\vec{q}\,)](q).$$
We write the explicit series expression using the fact that $Q_l(\vec{q}\,)=c_l\vec{q}\,^l$ for all $ l\geq 0$ and obtain
 \[ \begin{split}
 \displaystyle (Q_k\odot Q_s)(q)& =\sum _{j=0}^\infty \frac{(-1)^jx_0^j}{j!}\partial_{\vec{q}\,}^j \left( Q_k(\vec{q}\,)Q_s(\vec{q}\,) \right) \\
&= c_kc_s\sum _{j=0}^\infty \frac{(-1)^jx_0^j}{j!}\partial_{\vec{q}\,}^j \left(\vec{q}\,^{k+s} \right).
 \\
\end{split}
\]
In particular, we get \begin{equation}\label{ckcs}
\displaystyle (Q_k\odot Q_s)(q)=c_kc_s CK\left(\vec{q}\,^{k+s}\right)(q), \textbf{ } q\in\Hq,  k,s\geq 0,
\end{equation}
with $c_l:=\displaystyle \sum_{j=0}^l(-1)^jT_j^l, \textbf{  } \forall l\geq 0.$
On the other hand, we observe that $Q_{k+s}$ is also Fueter regular on $\Hq$. Moreover, it is restriction to $x_0=0$ gives $$Q_{k+s}(\vec{q}\,)=c_{k+s} \vec{q}\,^{k+s}.$$
Therefore, by uniqueness of the C-K extension we get
\begin{equation}\label{CKk+s}
Q_{k+s}(q)=c_{k+s}CK\left(\vec{q}\,^{k+s}\right)(q),\textbf{ } \forall q\in\Hq.
\end{equation}
Hence, we combine \eqref{ckcs} and \eqref{CKk+s} to conclude that
$$\displaystyle (Q_k\odot Q_s)(q)=\frac{c_kc_s}{c_{k+s}}Q_{k+s}(q), \textbf{} \forall q\in\Hq, \forall k,s\geq 0.$$
\end{proof}
\begin{rem}
If we consider the Fueter regular polynomials given by $P_k=\dfrac{Q_k}{c_k}$, $\forall k\geq 0$. Then, the classical multiplication rule holds, in the sense that we have
\begin{equation}
P_k\odot P_s=P_{k+s}, \textbf{ } \forall k,s\geq 0.
\end{equation}
\end{rem}
\begin{cor}
Let $k,s \geq 0$. Then, for any $q=x_0+\vec{q}\,\in\Hq$ we have

$$\displaystyle (Q_k\odot Q_s)(q)=c_kc_s \lambda_{0}^{k+s}r^{k+s}\left(C_{k+s}^1\left(\frac{x_0}{r}\right)+\frac{2}{k+s+2}C_{k+s-1}^2\left(\frac{x_0}{r} \right)\frac{\vec{q}\,}{r}\right), $$
where $C_t^\nu$ are the Gegenbauer polynomials, $\lambda_0$ is a constant and $r^2=|q|^2$.
\end{cor}
\begin{proof}
Proposition \ref{QkCKQs} gives
 $$(Q_k\odot Q_s)(q\,)|_{x_0=0}=c_kc_s \vec{q}\,^{k+s}, \textbf{ }  k,s\geq 0,$$
thus, by the regularity of the C-K product $Q_k\odot Q_s$ and uniqueness of the C-K extension we have that$$(Q_k\odot Q_s)(q)=c_kc_s CK[\vec{q}\,^{k+s}], \textbf{ }q\in\Hq, k,s\geq 0.$$
Hence, the result follows as a direct application of Theorem 2.2.1 in \cite{DSS1992} that gives the expression of the C-K extension for the vector part powers in terms of Gegenbauer polynomials.
\end{proof}
\begin{rem}
We note that the Appell polynomials given by \eqref{Qk2} define a family of Fueter regular functions of axial type (or axially Fueter regular functions), in the sense that if we write $q=x_0+\omega|\vec{q}\,|\in \Omega$ with $\omega\in\Sq$ there exist two quaternionic valued functions $A=A(x_0,|\vec{q}\,|)$ and $B=B(x_0,|\vec{q}\,|)$ independent of $\omega$ such that we have
\begin{equation}
Q_k(q)=A(x_0,|\vec{q}\,|)+\omega B(x_0,|\vec{q}\,|), \textbf{ } \forall k\geq 0.
\end{equation}
\end{rem}

We end this section by proving a converse result of the previous remark. This allows to characterize axially Fueter regular functions on quaternionic axially symmetric slice domains in terms of the Appell system $(Q_k)_{k\geq 0}$.
\begin{thm}
Let $\Omega\subseteq\Hq$ be an axially symmetric slice domain. Let $g$ be an axially Fueter regular function on $\Omega$. Then, there exist some quaternion coefficients $(\alpha_k)_{k\geq 0}$ such that we have the expansion

\begin{equation}
g(q)=\displaystyle \sum_{k=0}^\infty Q_k(q)\alpha_k, \textbf{ } \forall q\in\Omega.
\end{equation}
\end{thm}
\begin{proof}
We note that $g$ is an axially Fueter regular function on $\Omega$. Thus, by the inverse Fueter mapping theorem proved in \cite{ColomboSabadiniSommen2011} there will exist $f\in\mathcal{SR}(\Omega)$ such that we have
\begin{equation}\label{tauf1}
g=\tau(f),
\end{equation}
where $\tau=\Delta_{\R^4}$ is the Fueter mapping. Then, using the series expansion theorem for slice hyperholomorphic functions there exist some quaternion coefficients $(a_k)_{k\geq 0}$ so that we can write
\begin{equation} \label{se}
f(q)=\displaystyle \sum_{k=0}^\infty q^k a_k, \textbf{ } \forall q\in \Omega.
\end{equation}

In particular, we apply the Fueter mapping $\tau$ on \eqref{se} and get
$$\tau(f)(q)=\displaystyle \sum_{k=0}^\infty \tau(q^k)a_k.$$
However, we know by \cite{DKS2019} that $$\tau(q^k)=-2(k-1)kQ_{k-2}, \forall k\geq 2.$$
Therefore, we continue the calculations and obtain
\begin{equation}\label{tauf2}
\tau(f)=\displaystyle \sum_{k=0}^\infty Q_k\alpha_k,
\end{equation}
where we have set $\alpha_k=-2(k+1)(k+2)a_{k+2}, \textbf{ } \forall k\geq 0.$ Hence, comparing \eqref{tauf1} with \eqref{tauf2} we conclude that  $$g(q)=\displaystyle \sum_{k=0}^\infty Q_k(q)\alpha_k, \textbf{ } \forall q\in\Omega.$$
This ends the proof.
\end{proof}

\section{The Fock space case}

In this section, we consider the Clifford-Appell Fock space in the setting of quaternions which is defined by
$$\mathcal{F}(\Hq):=\left\lbrace{\displaystyle \sum_{k=0}^\infty Q_k\alpha_k; \textit{ } \alpha_k\in \Hq \ :\  \sum_{k=0}^\infty k!|\alpha_k|^2<\infty}\right\rbrace.$$

This space corresponds to the space $\mathcal{HM}_b$ in Definition \ref{DefHb} associated with the sequence $b=k!$, $k\geq 0$ on the domain $\Omega=\Hq$. Let $f=\displaystyle \sum_{k=0}^\infty Q_k\alpha_k \text{ and } g=\displaystyle \sum_{k=0}^\infty Q_k\beta_k$ in $\mathcal{F}(\Hq)$ we can equip $\mathcal{F}(\Hq)$ with the scalar product $$\scal{f,g}_{\mathcal{F}(\Hq)}=\displaystyle\sum_{k=0}^\infty k!\overline{\alpha_k}\beta_k.$$
Then, we can see that all the evaluation mappings on $\mathcal{F}(\Hq)$ are continuous. Indeed, we prove the following estimate
\begin{prop}\label{FockEst}
For any $f\in\mathcal{F}(\Hq)$ and $q\in\Hq$, we have
\begin{equation}
|f(q)|\leq e^{\frac{|q|^2}{2}}\Vert f \Vert_{\mathcal{F}(\Hq)}.
\end{equation}

\end{prop}
\begin{proof}
We write $f(q)=\displaystyle\sum_{k=0}^\infty Q_k(q)\alpha_k$. Thus, we have $$|f(q)|\leq \displaystyle \sum _{k=0}^\infty\frac{|Q_k(q)|}{\sqrt{k!}}|\alpha_k|\sqrt{k!}.$$
Then, by the Cauchy-Schwarz inequality we obtain $$|f(q)|\leq \left(\displaystyle \sum_{k=0}^\infty \frac{|Q_k(q)|^2}{k!}\right)^{\frac{1}{2}} \left(\displaystyle \sum_{k=0}^\infty k!|\alpha_k|^2\right)^{\frac{1}{2}}$$
However, we know that $|Q_k(q)|\leq |q|^k$ for all $q\in\Hq$ (see the proof of Proposition 4.5 in \cite{DKS2019}).  Hence, we get $$|f(q)|\leq e^{\frac{|q|^2}{2}}\Vert f \Vert_{\mathcal{F}(\Hq)}.$$

\end{proof}
As a consequence, we have the following result
\begin{thm}
The set $\mathcal{F}(\Hq)$ is a right quaternionic Hilbert space of Cauchy-Fueter regular functions whose reproducing kernel is given by $$K_{\mathcal{F}(\Hq)}(q,p)=\displaystyle \sum_{k=0}^\infty \frac{Q_k(q)\overline{Q_k(p)}}{k!}, \textbf{  } \forall (q,p)\in\Hq\times \Hq.$$
Moreover, if we set $\psi_k(q)=\displaystyle\frac{Q_k(q)}{\sqrt{k!}}, \textbf{  }  k\geq 0,$ then, the family $\lbrace{\psi_k}\rbrace_{k\geq 0} $ form an orthonormal basis of $\mathcal{F}(\Hq)$.
\end{thm}
\begin{proof}
For a fixed $p\in\Hq$, we consider the function defined by $$K_p(q)=\displaystyle \sum_{k=0}^\infty Q_k(q)\beta_k(p) , \text{ }  \forall q\in\Hq, \text{ where } \beta_k(p)=\displaystyle \frac{\overline{Q_k(p)}}{k!}.$$
We observe that $$\displaystyle \sum_{k=0}^\infty k!|\beta_k(p)|^2=\sum_{k=0}^\infty \frac{|Q_k(p)|^2}{k!}\leq e^{|q|^2}<\infty. $$
So, the function $K_p$ belongs to $\mathcal{F}(\Hq)$ for all $p\in\Hq$. Now, let $f= \displaystyle \sum_{k=0}^\infty Q_k\alpha_k $ be any function in $\mathcal{F}(\Hq)$. Then $$\scal{K_p,f}_{\mathcal{F}(\Hq)}=\displaystyle \sum_{k=0}^\infty k! \overline{\beta_k(p)}\alpha_k= \sum_{k=0}^\infty Q_k(p)\alpha_k=f(p), \textbf{  } \forall p\in\Hq,$$
therefore, the reproducing kernel of the space $\mathcal{F}(\Hq)$ is given by $$K_{\mathcal{F}(\Hq)}(q,p)=\displaystyle \sum_{k=0}^\infty \frac{Q_k(q)\overline{Q_k(p)}}{k!}, \textbf{  } \forall (q,p)\in\Hq\times \Hq.$$
It is clear by definition of the scalar product that $$\scal{\psi_k,\psi_j}_{\mathcal{F}(\Hq)}=\delta_{k,j} ,\textbf{  } \forall k,j\in \N.$$ Furthermore, let $f= \displaystyle \sum_{k=0}^\infty Q_k\alpha_k $ in $\mathcal{F}(\Hq)$ be such that $$\scal{\psi_k,f}_{\mathcal{F}(\Hq)}=0, \textbf{  } \forall k\in \N.$$
We have $$\displaystyle \sqrt{k!} \alpha_k=\scal{\psi_k,f}_{\mathcal{F}(\Hq)}=0, \textbf{  } \forall k\in \N,$$ so, $f=0 \text{ for any  } q\in\Hq.$ In particular, this proves that $\lbrace{\psi_k}\rbrace_{k\geq 0} $ form an orthonormal basis of $\mathcal{F}(\Hq)$.
\end{proof}
\begin{rem}
We note that  \begin{itemize}

   \item[i)]$K_{\mathcal{F}(\mathbb{H})}(\vec{q}\,,\vec{p})=\displaystyle \sum_{k=0}^\infty (-1)^k\frac{c_k^2}{k!}\vec{q}\,^k\vec{p}^k, \textbf{  } \forall (q,p)\in\mathbb{H}_0\times \mathbb{H}_0.$
   \item[ii)]$K_{\mathcal{F}(\mathbb{H})}(x,y)=\displaystyle e^{xy}, \textbf{  } \forall (x,y)\in \R\times \R.$

\end{itemize}
\end{rem}
Now we turn our attention to the notion of creation operator associated with the Clifford-Appell Fock space $\mathcal{F}(\Hq)$. For this, we consider a sequence of real numbers $\gamma=(\gamma_k)_{k\geq 0}$ that allows to define a weighted shift operator by \begin{equation}\label{shiftgamma}
T_\gamma(Q_k):=\gamma_k Q_{k+1}, \textbf{ } \forall k\geq 0.
\end{equation}
We would like to preserve in this setting the main properties of adjoint and commutation rules satisfied by the standard creation and annihilation operators on the Fock space.
First, we deal with the following
\begin{prop}
Let $\gamma$ be a sequence with $\gamma_0=1$ and such that \eqref{shiftgamma} is well defined. Then, we have $$\left[\displaystyle\frac{\overline{\partial}}{2} T_\gamma , T_\gamma\displaystyle\frac{\overline{\partial}}{2} \right]=\mathcal{I}_{\mathcal{F}(\Hq)},$$
if and only if $$\displaystyle\gamma_k=\frac{1+k\gamma_{k-1}}{1+k}, \textbf{  } \forall k\geq 1.$$
\end{prop}
\begin{proof}
Let $\displaystyle f=\sum_{k=0}^\infty Q_k\alpha_k$ be a function in $\mathcal{F}(\Hq)$. Then, we have
$$\displaystyle T_\gamma(f)=\sum_{k=0}^\infty\gamma_kQ_{k+1}\alpha_k \text{ and } \frac{\overline{\partial}}{2}(f)=\sum_{k=1}^\infty kQ_{k-1}\alpha_k.$$
Thus, we obtain
$$\displaystyle \frac{\overline{\partial}}{2}T_\gamma(f)=\sum_{k=0}^\infty(k+1)\gamma_kQ_{k}\alpha_k \text{ and } T_\gamma \frac{\overline{\partial}}{2}(f)=\sum_{k=1}^\infty k\gamma_{k-1} Q_{k}\alpha_k.$$
Therefore, it follows that
\begin{equation}\label{Comm}
\left[\displaystyle\frac{\overline{\partial}}{2} T_\gamma , T_\gamma\displaystyle\frac{\overline{\partial}}{2} \right](f)=\gamma_0Q_0\alpha_0+\sum_{k=1}^\infty [(k+1)\gamma_k-k\gamma_{k-1}]Q_k\alpha_k
\end{equation}
We can see that if  $$\displaystyle\gamma_k=\frac{1+k\gamma_{k-1}}{1+k}, \textbf{  } \forall k\geq 1,$$
we have then $$(k+1)\gamma_k-k\gamma_{k-1}=1, \textbf{ } \forall k\geq 1.$$
Therefore, using the condition $\gamma_0=1$ and formula \eqref{Comm} we obtain $$\left[\displaystyle\frac{\overline{\partial}}{2} T_\gamma , T_\gamma\displaystyle\frac{\overline{\partial}}{2} \right](f)=Q_0\alpha_0+\sum_{k=1}^\infty Q_k\alpha_k=f.$$
For the converse, if we assume that  $$\left[\displaystyle\frac{\overline{\partial}}{2} T_\gamma , T_\gamma\displaystyle\frac{\overline{\partial}}{2} \right](f)=f,$$ we apply \eqref{Comm} and get
$$\gamma_0Q_0(q)\alpha_0+\sum_{k=1}^\infty [(k+1)\gamma_k-k\gamma_{k-1}]Q_k(q)\alpha_k=\sum_{k=0}^\infty Q_k(q)\alpha_k, \textbf{  } \forall q\in\Hq.$$

In particular, using the fact that $Q_k(t)=t^k, \forall t\in\R$ and $\gamma_0=1$ we observe that
$$\alpha_0+\sum_{k=1}^\infty [(k+1)\gamma_k-k\gamma_{k-1}]t^k\alpha_k=\sum_{k=0}^\infty t^k\alpha_k, \textbf{  } \forall t\in\R.$$
Therefore, comparing the coefficients of the same degree we obtain $$(k+1)\gamma_k-k\gamma_{k-1}=1, \textbf{ } \forall k\geq 1.$$
Hence, we have the condition   $$\displaystyle\gamma_k=\frac{1+k\gamma_{k-1}}{1+k}, \textbf{  } \forall k\geq 1.$$
\end{proof}
Furthermore, we can prove the following
\begin{prop}\label{gamma1}
Let $\gamma$ be a sequence with $\gamma_0=1$ and such that \eqref{shiftgamma} holds.  If one of the following properties is satisfied
\begin{enumerate}
\item[i)] $\left[\displaystyle\frac{\overline{\partial}}{2} T_\gamma , T_\gamma\displaystyle\frac{\overline{\partial}}{2} \right]=\mathcal{I}_{\mathcal{F}(\Hq)};$
\item[ii)]  $T_\gamma$ is the adjoint of the hypercomplex derivative $\displaystyle\frac{\overline{\partial}}{2}$;
\end{enumerate}

then, we have $$\displaystyle\gamma_k=1, \textbf{  } \forall k\geq 0.$$
\end{prop}
\begin{proof}
We observe that condition  i) and Proposition \ref{gamma1} show that
$$\displaystyle\gamma_k=\frac{1+k\gamma_{k-1}}{1+k}, \textbf{  } \forall k\geq 1.$$Thus, since $\gamma_0=1$ a simple induction reasoning allows to prove that if i) holds then $\gamma_k=1$, for all $k\geq 1$.
On the other hand, the condition ii) implies in particular that we have $$\scal{\frac{\overline{\partial}}{2}(Q_k),Q_j}_{\mathcal{F}(\Hq)}=\scal{Q_k,T_\gamma(Q_j)}_{\mathcal{F}(\Hq)}, \textbf{ } \forall k,j\geq 1.$$
So, we conclude $$k(k-1)!\delta_{k-1,j}=\gamma_j k!\delta_{k,j+1}, \textbf{  } \forall k,j\geq 1,$$
where $\delta_{m,n}$ is the Kronecker symbol.
In particular, this leads to the same conclusion that $\gamma_j=1$, $j\geq 1$.
\end{proof}

\begin{rem}
We note that thanks to Proposition \ref{gamma1} the only operator $T_\gamma$ that can play the role of the creation operator with respect to the Clifford-Appell system should act as follows
\begin{equation}\label{creapp}
T_\gamma(Q_k)=Q_{k+1}, \textbf{ } \forall k\geq 0.
\end{equation}
\end{rem}
We now introduce the notion of creation operator associated with the quaternionic Hilbert space $\mathcal{HM}_b$ in terms of the C-K product that allows to have the property \eqref{creapp} . To this end, let $k\geq 0$,  and we define first the family of operators given by \begin{equation}
\displaystyle \mathcal{S}_k(f):=\frac{c_{1+k}}{c_1c_k}Q_1\odot f, \textbf{ }\forall f\in \mathcal{HM}_b
\end{equation}
where $\odot$ denote the C-K product and $c_l:=\displaystyle \sum_{j=0}^l(-1)^jT_j^l, \textbf{  } \forall l\geq 0.$ \\
Then, for $f=\displaystyle \sum_{k=0}^\infty Q_k\alpha_k $ in $\mathcal{HM}_b$ we consider the operator $\mathcal S$ defined by applying $\mathcal{S}_k$ on each component with the corresponding degree as follows

\begin{equation}\label{Screation}
\mathcal{S}(f):=\displaystyle \sum_{k=0}^\infty \mathcal{S}_k(Q_k)\alpha_k.
\end{equation}
Therefore, we have the explicit expression given by
\begin{equation}\label{Sexp}
\mathcal{S}(f):=\displaystyle \frac{1}{c_1}\sum_{k=0}^\infty \frac{c_{1+k}}{c_k}[Q_1\odot Q_k]\alpha_k.
\end{equation}
We note that the operator $\mathcal{S}$ acts like the classical shift operator with respect to the Clifford-Appell system $(Q_k)_{k\geq 0}$. This can be seen in the following
\begin{prop}\label{shift}
For all $k\geq 0$, we have $$\mathcal{S}(Q_k)(q)=Q_{k+1}(q), \textbf{ } \forall q\in \Hq.$$
\end{prop}
\begin{proof}
Let $k\geq 0$. Then, for all $q\in\Hq$ we have
 \[ \begin{split}
 \displaystyle \mathcal{S}(Q_k)(q)& =\mathcal{S}_k(Q_k)(q) \\
&= \frac{c_{1+k}}{c_1c_k}(Q_1\odot Q_k)(q).
 \\
\end{split}
\]
Now, we apply Proposition \ref{QkCKQs} and get $$Q_1\odot Q_k=\frac{c_1c_k}{c_{1+k}}Q_{k+1}.$$
Hence, we obtain $$\mathcal{S}(Q_k)=Q_{k+1}.$$
\end{proof}
As a consequence of Proposition \ref{shift} we note that the creation operator on $\mathcal{F}(\Hq)$ given by \eqref{Screation} acts as follows $$\mathcal{S}(\displaystyle \sum_{k=0}^\infty Q_k\alpha_k)=\displaystyle \sum_{k=0}^\infty Q_{k+1}\alpha_k.$$
The annihilation operator corresponds to  the hypercomplex derivative $$\frac{\overline{\partial}}{2}:=\displaystyle\frac{1}{2}\left(\frac{\partial}{\partial x_0}-i\frac{\partial}{\partial x_1}-j\frac{\partial}{\partial x_2}-k\frac{\partial}{\partial x_3}\right).$$
It is known by the Appell property that $$\frac{\overline{\partial}}{2}(Q_k)=kQ_{k-1}, \textbf{  }\forall k\geq 1.$$ The domains of $\mathcal{S}$ and $\displaystyle\frac{\overline{\partial}}{2}$ in $\mathcal{F}(\Hq)$ are denoted respectively by
$$D(\mathcal{S}):=\lbrace{f\in\mathcal{F}(\Hq); \textbf{  } \mathcal{S}(f)\in \mathcal{F}(\Hq)}\rbrace $$  and $$ D(\frac{\overline{\partial}}{2}):=\lbrace{f\in\mathcal{F}(\Hq); \textbf{  } \frac{\overline{\partial}}{2}(f)\in \mathcal{F}(\Hq)}\rbrace.$$
We note that the creation operator $\mathcal{S}$ and the hypercomplex derivative $\frac{\overline{\partial}}{2}$  are quaternionic right linear operators densely defined on $\mathcal{F}(\Hq)$ since $\left\lbrace{\frac{Q_k}{\sqrt{k!}}}\right\rbrace_{k\geq 0}$ is an orthonormal basis of the quaternionic Fock Hilbert space.  In the sequel, we shall prove some different properties of these operators:
\begin{prop}
$\mathcal{S}$ and  $\frac{\overline{\partial}}{2}$ are two closed quaternionic operators on $\mathcal{F}(\Hq)$.
\end{prop}
\begin{proof}
We consider the graph of $\mathcal{S}$ defined by $$\mathcal{G}(\mathcal{S}):=\lbrace{(f,\mathcal{S} f); f\in D(\mathcal{S})}\rbrace.$$  Let us show that $\mathcal{G}(\mathcal{S})$ is closed. Indeed, let $\phi_n$ be a sequence in $D(\mathcal{S})$ such that $\phi_n$ and $\mathcal{S}\phi_n$ converge to $\phi$ and $\psi$ respectively on $\mathcal{F}(\Hq)$. Then, thanks to Proposition \ref{FockEst} we have $$\vert{\phi_n(q)-\phi(q)}\vert\leq e^{\frac{|q|^2}{2}} \Vert{\phi_n-\phi}\Vert_{\mathcal{F}(\Hq)}$$
and
 $$
 \vert{\mathcal{S}\phi_n(q)-\psi(q)}\vert\leq e^{\frac{|q|^2}{2}} \Vert{\mathcal{S}\phi_n-\psi}\Vert_{\mathcal{F}(\Hq)}. $$
Therefore, it follows that $\phi_n$ and $\mathcal{S}\phi_n$ converge  pointwise to $\phi$ and $\psi$, respectively. This leads to $\psi=\mathcal{S}\phi$ which ends the proof. The same technique could be adapted to prove the closedness of the hypercomplex derivative on $\mathcal{F}(\Hq)$.
\end{proof}
Furthermore, we prove also the following

\begin{prop}
Let $f\in\mathcal{F}(\Hq)$. Then, $\mathcal{S} (f)$ belongs to $\mathcal{F}(\Hq)$ if and only if $\displaystyle\frac{\overline{\partial}}{2} f$ belongs to $\mathcal{F}(\Hq)$. In particular,  this means that we have  $$D(\mathcal{S})=D\left(\displaystyle\frac{\overline{\partial}}{2}\right).$$
\end{prop}
\begin{proof}
We write $f=\displaystyle \sum_{k=0}^\infty Q_k\alpha_k$ in $\mathcal{F}(\Hq)$. Then, we have $$\mathcal{S}(f)=\displaystyle \sum_{h=1}^\infty Q_h\alpha_{h-1}. $$
In particular, we have
\begin{equation}\label{Sfnorm}
 ||\mathcal{S}(f)||_{\mathcal{F}(\Hq)}^2=\sum_{h=1}^\infty h!|\alpha_{h-1}|^2.
\end{equation}

On the other hand, using the Appell property with respect to the hypercomplex derivative we have
$$\displaystyle \frac{\overline{\partial}}{2}(f)=\sum_{h=0}^\infty Q_h\beta_h, \textbf{ } \beta_h=(h+1)\alpha_{h+1}, \forall h\geq 0.$$

Some calculations lead to \begin{equation}\label{hydernorm}
||\displaystyle \frac{\overline{\partial}}{2}(f)||_{\mathcal{F}(\Hq)}^2=\sum_{h=1}^\infty h (h!)|\alpha_{h}|^2 .
\end{equation}
We note that by \eqref{Sfnorm} we have
\[ \begin{split}
 \displaystyle ||\mathcal{S} f||_{\mathcal{F}(\Hq)}^2 & =\sum _{h=0}^\infty (h+1)!|\alpha_h|^2\\
&= \sum _{h=0}^\infty (h+1)h!|\alpha_h|^2\\ &=\sum _{h=0}^\infty h(h)!|\alpha_h|^2+\sum _{h=0}^\infty h!|\alpha_h|^2 .
\end{split}
\]

Therefore, we use \eqref{hydernorm}  in order to get
\begin{equation}\label{^}
 \displaystyle ||\mathcal{S} f||_{\mathcal{F}(\Hq)}^2=||\frac{\overline{\partial}}{2} f||_{\mathcal{F}(\Hq)}^2+||f||_{\mathcal{F}(\Hq)}^2.
\end{equation}

Hence, formula \eqref{^} shows that $||\mathcal{S} f||_{\mathcal{F}(\Hq)}<\infty$ if and only if $||\frac{\overline{\partial}}{2} f||_{\mathcal{F}(\Hq)}<\infty$ which ends the proof.
\end{proof}
Now, we prove the adjoint property
\begin{prop}
Let $f\in D(\displaystyle\frac{\overline{\partial}}{2})$ and $g\in D(\mathcal{S})$. Then, we have
$$\scal{\displaystyle\frac{\overline{\partial}}{2}f,g}_{\mathcal{F}(\Hq)}=\scal{f,\mathcal{S}(g)}_{\mathcal{F}(\Hq)}.$$

\end{prop}
\begin{proof}
Let $\displaystyle f=\sum _{k=0}^\infty Q_k\alpha_k$ in $D(\displaystyle\frac{\overline{\partial}}{2})$ and $\displaystyle g=\sum _{k=0}^\infty Q_k\beta_k$ in  $D(\mathcal{S})$. Thus, we have  \[ \begin{split}
 \displaystyle \frac{\overline{\partial}}{2}f & =\sum _{k=0}^\infty \frac{\overline{\partial}}{2}(Q_k)\alpha_k \\
&= \sum _{k=1}^\infty kQ_{k-1} \alpha_k\\ &=\sum _{h=0}^\infty (h+1)Q_{h} \alpha_{h+1}.
\end{split}
\]
On the other hand, making use of Proposition \ref{shift} we have
 \[ \begin{split}
 \displaystyle \mathcal{S}(g) & =\sum _{k=0}^\infty \mathcal{S}(Q_k)\beta_k \\
&= \sum _{k=0}^\infty Q_{k+1} \beta_k\\ &=\sum _{k=1}^\infty Q_{k} \beta_{k-1}.
\end{split}
\]
Therefore, we obtain  $$\scal{\displaystyle\frac{\overline{\partial}}{2}f,g}_{\mathcal{F}(\Hq)}=\sum_{k=0}^\infty (k+1)!\overline{\alpha_{k+1}}\beta_k=\scal{f,\mathcal{S}(g)}_{\mathcal{F}(\Hq)}.$$
This ends the proof.
\end{proof}
\begin{prop}
Let $f\in \mathcal{D}(\displaystyle\frac{\overline{\partial}}{2})\cap\mathcal{D}(\mathcal{S}) $. Then, we have $$\left[\displaystyle\frac{\overline{\partial}}{2} \mathcal{S} , \mathcal{S}\displaystyle\frac{\overline{\partial}}{2} \right](f)=f.$$

\end{prop}
\begin{proof}
Let $\displaystyle f=\sum _{k=0}^\infty Q_k\alpha_k$ be in $ \mathcal{D}(\displaystyle\frac{\overline{\partial}}{2})\cap\mathcal{D}(\mathcal{S})$. Thus, computations using Proposition \ref{shift} and the Appell property give $$\displaystyle\frac{\overline{\partial}}{2}\mathcal{S}(f)=\sum_{k=0}^\infty (k+1)Q_k\alpha_k\ \ \text{ and }\ \  \mathcal{S}\frac{\overline{\partial}}{2}(f)=\sum_{k=0}^\infty kQ_k\alpha_k.$$ In particular, it shows that $$\displaystyle\frac{\overline{\partial}}{2}\mathcal{S}(f)-\mathcal{S}\frac{\overline{\partial}}{2}(f)=f.$$ This ends the proof.
\end{proof}
\begin{rem}
Note that the creation and annihilation operators denoted respectievly by $\mathcal{S}$ and $\displaystyle\frac{\overline{\partial}}{2} $ are adjoint of each other and satisfy the classical commutation rules on the Fock space of Fueter regular functions $\mathcal{F}(\Hq)$ like in the classical complex case. Moreover,  observe that we have also $\mathcal{S}\displaystyle\frac{\overline{\partial}}{2}(Q_k)=kQ_k$, for any $k\geq 1$. This property is related to the notion of number operators that appears in quantum mechanics.
\end{rem}
Let $(\eta_n)_{n\in\N}$ denote the normalized Hermite functions.
In order to study the Segal-Bargmann transform notion in this framework we introduce the Fueter regular kernel function given by
\begin{equation}
\mathcal{A}_\Hq^F(q,x):=\displaystyle \sum_{k=0}^\infty\frac{Q_k(q)}{\sqrt{k!}}\eta_k(x), \textbf{ } \forall (q,x)\in \Hq\times \R.
\end{equation}
Then, for any quaternionic valued function $\varphi$ in $L^2(\R,\Hq)$ and $q\in\Hq$ we define
\begin{equation}
\displaystyle \mathcal{B}_\Hq^F(\varphi)(q)=\int_\R \mathcal{A}_\Hq^F(q,x) \varphi(x)dx.
\end{equation}
We shall prove the following result:
\begin{thm} \label{SBT}
The integral transform $\mathcal{B}_\Hq^F$ defines an isometric isomorphism mapping the standard Hilbert space $L^2(\R,\Hq)$ onto the Clifford-Appell Fock space $ \mathcal{F}(\Hq)$.
\end{thm}
\begin{proof}
Let $\varphi \in L^2(\R,\Hq)$. We write $\varphi=\displaystyle \sum_{j=0}^\infty \eta_j(x)\beta_j$ such that $\Vert{\varphi}\Vert^2_{L^2(\R,\Hq)}=\displaystyle \sum_{j=0}^\infty |\beta_j|^2<\infty.$  Then, note that we have $$\displaystyle \mathcal{B}_\Hq^F(\varphi)(q)= \sum_{k=0}^\infty\frac{Q_k(q)}{\sqrt{k!}}\int_\R\eta_k(x)\varphi(x)dx.$$ So, by setting $\alpha_k=\displaystyle \frac{1}{\sqrt{k!}}\int_\R\eta_k(x)\varphi(x)dx$ for all $k\geq 0$, we get
 \[ \begin{split}
 \displaystyle \Vert \mathcal{B}_\Hq^F(\varphi) \Vert_{ \mathcal{F}(\Hq)}^2 & =\sum _{k=0}^\infty k!|\alpha_k|^2 \\
&= \sum _{k=0}^\infty \left|\int_\R \eta_k(x)\varphi(x)dx\right|^2.
\end{split}
\]
However, by definition of $\varphi$ and using the orthogonality of Hermite functions we obtain
 \[ \begin{split}
 \displaystyle \int_\R \eta_k(x)\varphi(x)dx  & =\sum _{j=0}^\infty \beta_j\int_\R\eta_k(x)\eta_j(x)dx =\beta_k, \textbf{ } \forall k\geq 0.
\end{split}
\]
Hence, we conclude that $$\displaystyle \Vert \mathcal{B}_\Hq^F(\varphi) \Vert_{ \mathcal{F}(\Hq)}^2=\sum_{j=0}^\infty |\beta_j|^2=\Vert{\varphi}\Vert^2_{L^2(\R,\Hq)}.$$ Moreover, observe that $$\displaystyle \mathcal{B}_\Hq^F(\eta_k)=\frac{Q_k}{\sqrt{k!}}, \text{  } \forall k\geq 0.$$
In particular, this allows to prove that  $\mathcal{B}_\Hq^F$ is an isometric isomorphism mapping the standard Hilbert space $L^2(\R,\Hq)$ onto the Fock space $ \mathcal{F}(\Hq)$ on the quaternions.
\end{proof}
Now, we consider the following:
\begin{pb}

 Is it possible to map $\mathcal{F}_{Slice}(\Hq)$ onto $\mathcal{F}(\Hq)$ without using the Fueter mapping, see \cite{DKS2019}, and keeping the isometry property ?

\end{pb}
To answer the question, we will compute $\mathcal{B}_\Hq^F$ composed with the slice hyperholomorphic Segal-Bargmann transform.

In order to answer this problem, we need the slice hyperholomorphic Segal-Bargmann transform given by \eqref{BargTra}.

Notice that thanks to these integral transforms $\mathcal{B}_\Hq^S$ and $\mathcal{B}_\Hq^F$ it is possible to relate the two notions of Fock spaces on the quaternions, namely the slice hyperholomorphic $ \mathcal{F}_{Slice}(\Hq)$ and the Cauchy-Fueter regular one $\mathcal{F}(\Hq).$ Indeed, for a fixed $i\in\Sq$, $f\in \mathcal{F}_{Slice}(\Hq)$ and $q\in\Hq$ we define the integral transform given by $$\Upsilon(f)(q):=\displaystyle \int_{\C_i}\mathcal{L}(q,z)f_i(z)d\mu_i(z),$$
where $\displaystyle d\mu_i(z):=\frac{1}{\pi}e^{-|z|^2}dA_i(z)$ and the kernel function is obtained by taking the series $$\displaystyle \mathcal{L}(q,z)=\sum_{k=0}^\infty\frac{Q_k(q)}{k!}\overline{z}^k, \textbf{  }\forall (q,z)\in\Hq\times \C_i.$$ Then, we prove:
\begin{thm}
The quaternionic integral transform $\Upsilon$ does not depend on the choice of the imaginary unit $i\in\mathbb S$. Furthermore, it defines an isometric isomorphism mapping the slice hyperholomrphic Fock space $\mathcal{F}_{Slice}(\Hq)$ onto the Clifford-Appell Fock space $\mathcal{F}(\Hq)$.

\end{thm}
\begin{proof}
Let $f\in\mathcal{F}_{Slice}(\Hq)$, by Proposition 3.11  in \cite{AlpayColomboSabadini2014} we have $$f(q)=\displaystyle \sum_{k=0}^\infty q^ka_k \text{ and } \sum_{k=0}^\infty |a_k|^2k!< \infty.$$
In particular, by definition of $\Upsilon$ we have
 \[ \begin{split}
 \displaystyle \Upsilon(f)(q) & =\int_{\C_i}\left(\sum _{k=0}^\infty \frac{Q_k(q)}{k!}\overline{z}^k\right)\left(\sum _{j=0}^\infty z^ja_j \right) d\mu_i(z)\\
&= \sum _{k,j=0}^\infty \frac{Q_k(q)}{k!}\left( \int_{\C_i} \overline{z}^kz^jd\mu_i(z)\right) a_j.
\end{split}
\]
However, it is known that $$\displaystyle  \int_{\C_i} \overline{z}^kz^jd\mu_i(z)=k!\delta_{k,j}.$$
Therefore, we get $$\Upsilon(f)(q)=\displaystyle \sum_{k=0}^\infty Q_k(q)a_k.$$
Hence, since the coefficients $(a_k)_{k\geq 0}$ do not depend on the choice of the imaginary unit $i$ we conclude that $\Upsilon(f)$ is well defined and does not depend on the choice of the imaginary unit. Now, we observe that the operator $\Upsilon$ can be obtained thanks to the commutative diagram such that we have $$\Upsilon= \mathcal{B}_\Hq^F\circ (\mathcal{B}_\Hq^S)^{-1}.$$ Indeed, to prove this fact. Let $f\in\mathcal{F}_{Slice}(\Hq)$ and set $$\displaystyle \phi(x)=(\mathcal{B}_\Hq^S)^{-1}(f)(x)=\int_{\C_i}\mathcal{A}_\Hq^S(\overline{z},x)f_i(z)d\mu_i(z).$$
Thus, for any $q\in\Hq$ we have:
 \[ \begin{split}
 \displaystyle \mathcal{B}_\Hq^F(\phi)(q) & =\int_{\C_i} \mathcal{A}_\Hq^F(q,x)\phi(x)dx\\
&= \int_{\C_i} \mathcal{A}_\Hq^F(q,x)\left(\int_{\C_i}\mathcal{A}_\Hq^S(\overline{z},x)f_i(z)d\mu_i(z)\right)dx\\ &=\int_{\C_i}\left(\int_\R \mathcal{A}_\Hq^F(q,x)\mathcal{A}_\Hq^S(\overline{z},x)dx\right)f_i(z)d\mu_i(z).
\end{split}
\]
Then, we set $$H(q,z)=\int_\R \mathcal{A}_\Hq^F(q,x)\mathcal{A}_\Hq^S(\overline{z},x)dx, \textbf{  } \forall (q,z)\in\Hq\times \C_i.$$
 So, for all $(q,z)\in\Hq\times \C_i$ we have
 \[ \begin{split}
 \displaystyle H(q,z) & =\int_{\C_i}\left(\sum_{k=0}^\infty \frac{Q_k(q)}{\sqrt{k!}}\eta_k(x)\right)\left(\sum_{j=0}^\infty \frac{\overline{z^j}}{\sqrt{j!}}\eta_j(x)\right)dx
 \\
&=\sum_{k,j=0}^\infty \frac{Q_k(q)}{\sqrt{k!}}\left(\int_\R\eta_k(x)\eta_j(x)dx \right)\frac{\overline{z^j}}{\sqrt{j!}}
\end{split}
\]
Then, using the fact that Hermite functions form an orthonormal basis of $L^2(\R,\Hq)$ we get
$$\displaystyle H(q,z)=\sum_{k=0}^\infty\frac{Q_k(q)}{k!}\overline{z}^k=\mathcal{L}(q,z), \textbf{  } \forall (q,z)\in\Hq\times \C_i.$$
At this stage, we replace $H(q,z)$ by its expression and conclude that  we have $$\Upsilon= \mathcal{B}_\Hq^F\circ (\mathcal{B}_\Hq^S)^{-1}.$$ Therefore, since both of $\mathcal{B}_\Hq^F$ and $\mathcal{B}_\Hq^S$ are isometric isomorphisms mapping $L^2(\R,\Hq)$ respectievly onto $\mathcal{F}(\Hq)$ and $\mathcal{F}_{Slice}(\Hq)$. This ends the proof.
 \end{proof}
This quaternionic operator satisfies also the following properties :
\begin{prop}\label{Rrange}
For all $n\geq 0$, we set $f_n(q)=\displaystyle \frac{q^n}{\sqrt{n!}}$ and $\phi_n(q)=\displaystyle \frac{Q_n(q)}{\sqrt{n!}},$  $q\in \Hq.$ Then, we have
\begin{enumerate}
\item[i)] $\Upsilon(f_n)=\phi_n,$  $\forall n\geq 0$.
\item[ii)] $\displaystyle \int_{\C_i} \mathcal{L}(q,z)\overline{\mathcal{L}(p,z)}d\mu_i(z)=K_{\mathcal{F}(\Hq)}(q,p), \textbf{  } \forall (q,p)\in\Hq\times\Hq.$
\end{enumerate}
\end{prop}
\begin{proof}
The first statement is a direct consequence of the fact that $$\Upsilon= \mathcal{B}_\Hq^F\circ (\mathcal{B}_\Hq^S)^{-1}.$$ This is combined with the two following relations $$(\mathcal{B}_\Hq^S)^{-1}(\eta_n)=f_n \text{  and  } \mathcal{B}_\Hq^F(f_n)=\phi_n, \text{  } \forall n\geq 0.$$
Now, let $(q,p)\in\Hq\times\Hq$. Then, we have  \[ \begin{split}
 \displaystyle \displaystyle \int_{\C_i} \mathcal{L}(q,z)\overline{\mathcal{L}(p,z)}d\mu_i(z) & =\sum_{k,j=0}^\infty \frac{Q_k(q)}{k!}\left( \int_{\C_i} \overline{z}^kz^jd\mu_i(z)\right)\frac{\overline{Q_j(p)}}{j!}\\
&=\sum_{k=0}^\infty \frac{Q_k(q)\overline{Q_k(p)}}{k!}, \\ &=K_{\mathcal{F}(\Hq)}(q,p).
\end{split}
\]
\end{proof}
\begin{cor}
Let $i\in\Sq$. Then, for all $x,y\in\R$ and $n\geq 0$, we have the following identities

\begin{enumerate}
\item[i)] $\displaystyle \int_{\C_i}e^{x\overline{z}}z^nd\mu_i(z)=x^n$.
\item[ii)] $\displaystyle \int_{\C_i}e^{x\overline{z}+yz}d\mu_i(z)=e^{xy}$.
\end{enumerate}

\end{cor}
\begin{proof}
Observe that we have \begin{equation}\label{-}
\mathcal{L}(t,z)=e^{t\overline{z}}, \textbf{  } \forall (t,z)\in\R\times\C_i.
\end{equation}
The first identity follows from i) of Proposition \ref{Rrange} combined with \eqref{-}.

The second statement is also a consequence of  \eqref{-} combined with ii) of Proposition \ref{Rrange} and the fact that $$K_{\mathcal{F}(\Hq)}(x,y)=e^{xy}, \textbf{  } \forall (x,y)\in\R\times\R.$$

\end{proof}
\section{The Hardy space case}
In this section, we study on the quaternionic unit ball $\Omega=\mathbb{B}$ the spaces associated to some sequence $b$ as considered in Definition \ref{DefHb}. First, we give some general proofs related to these spaces $\mathcal{HM}_b(\mathbb{B}).$ Then, we will give more specific results on the Clifford-Appell Hardy space in this framework that corresponds to the sequence $b_k=1, \forall k\geq 0$.
In all this part, we take $\Omega=\mathbb{B}$ and $b=(b_k)_{k\geq 0}$  a non decreasing sequence with $b_0=1$. Then, we have
\begin{prop}
The following estimate holds
$$
  \displaystyle |f(q)|\leq \left(\sum_{k=0}^\infty \frac{|q|^{2k}}{b_k}\right)^{\frac{1}{2}} \Vert f \Vert_{\mathcal{HM}_b}, \textbf{ } f\in\mathcal{HM}_b(\mathbb{B}), \textbf{ }q\in\mathbb{B}.$$
\end{prop}
\begin{proof}
Let us consider $f(q)=\displaystyle\sum_{k=0}^\infty Q_k(q)\alpha_k$ in $\mathcal{HM}_b(\mathbb{B})$. Thus, we have $$|f(q)|\leq \displaystyle \sum _{k=0}^\infty\frac{|Q_k(q)|}{\sqrt{b_k}}|\alpha_k|\sqrt{b_k}.$$
Then, by the Cauchy-Schwarz inequality we have $$|f(q)|\leq \left(\displaystyle \sum_{k=0}^\infty \frac{|Q_k(q)|^2}{b_k}\right)^{\frac{1}{2}} \left(\displaystyle \sum_{k=0}^\infty b_k|\alpha_k|^2\right)^{\frac{1}{2}}$$
However, we know that $|Q_k(q)|\leq |q|^k$.  Hence, we get $$\displaystyle |f(q)|\leq \left(\sum_{k=0}^\infty \frac{|q|^{2k}}{b_k}\right)^{\frac{1}{2}} \Vert f \Vert_{\mathcal{HM}_b}.$$
\end{proof}
As a consequence, we get this result

\begin{thm}
The sets $\mathcal{HM}_b(\mathbb{B})$ are right quaternionic reproducing kernel Hilbert spaces.  Their reproducing kernel functions are given by
 \begin{equation}
  K_{\mathcal{H}_b(\mathbb{B})}(q,p)=\displaystyle \sum_{k=0}^\infty \frac{Q_k(q)\overline{Q_k(p)}}{b_k}, \textbf{  } \forall (q,p)\in\mathbb{B}\times \mathbb{B}.
 \end{equation}

Furthermore, the family  $\displaystyle\lbrace \psi_k^b:=\frac{Q_k}{\sqrt{b_k}}, \textbf{  }k\geq 0 \rbrace$ forms an orthonormal basis of $\mathcal{H}_b(\mathbb{B})$.
\end{thm}
\begin{proof}
For a fixed $p\in\mathbb{B}$, we consider the function defined by $$K_p(q)=\displaystyle \sum_{k=0}^\infty Q_k(q)\beta_k(p) , \text{ }  \forall q\in\mathbb{B}, \text{ where } \beta_k(p)=\displaystyle \frac{\overline{Q_k(p)}}{b_k}.$$
Thanks to  the d'Alembert ratio test for power series, we have $$\displaystyle \sum_{k=0}^\infty b_k|\beta_k(p)|^2=\sum_{k=0}^\infty \frac{|Q_k(p)|^2}{b_k}\leq \sum_{k=0}^\infty \frac{|q|^{2k}}{b_k} < \infty. $$
So, the function $K_p$ belongs to $\mathcal{H}_b(\mathbb{B})$ for any $p\in\mathbb{B}$. Now, let $f= \displaystyle \sum_{k=0}^\infty Q_k\alpha_k \in\mathcal{H}_b(\mathbb{B})$. Then, we have $$\scal{K_p,f}_{\mathcal{H}(\mathbb{B})}=\displaystyle \sum_{k=0}^\infty b_k \overline{\beta_k(p)}\alpha_k= \sum_{k=0}^\infty Q_k(p)\alpha_k=f(p), \textbf{  } \forall p\in\mathbb{B}.$$
Therefore, the reproducing kernel of the space $\mathcal{H}_b(\mathbb{B})$ is given by $$K_{\mathcal{H}_b(\mathbb{B})}(q,p)=\displaystyle \sum_{k=0}^\infty \frac{Q_k(q)\overline{Q_k(p)}}{b_k}, \textbf{  } \forall (q,p)\in\mathbb{B}\times \mathbb{B}.$$
It is clear by definition of the scalar product that $$\scal{\psi_k^b,\psi_j^b}_{\mathcal{H}_b(\mathbb{B})}=\delta_{k,j} ,\textbf{  } \forall k,j\in \N.$$ Furthermore, let $f= \displaystyle \sum_{k=0}^\infty Q_k\alpha_k $ in $\mathcal{H}_b(\mathbb{B})$ be such that $$\scal{\psi_k^b,f}_{\mathcal{H}_b(\mathbb{B})}=0, \textbf{  } \forall k\in \N.$$
Thus, we have $$\displaystyle \sqrt{b_k} \alpha_k=\scal{\psi_k^b,f}_{\mathcal{H}_b(\mathbb{B})}=0, \textbf{  } \forall k\in \N.$$ So, $f=0 \text{ for any  } q\in\mathbb{B}.$ In particular, this proves that $\lbrace{\psi_k^b}\rbrace_{k\geq 0} $ form an orthonormal basis of $\mathcal{HM}_b(\mathbb{B})$.
\end{proof}
\begin{rem}
The Clifford-Appell Hardy space corresponds to the sequence $b$ with all the terms equal to $1$,  and will be denoted  simply $\mathcal{H}(\mathbb{B})$. In this case, the previous results of this section read as follows \begin{itemize}
\item[i)] $
  \displaystyle |f(q)|\leq  \frac{\Vert f \Vert_{\mathcal{H}(\mathbb{B})}}{\left(1-|q|^2\right)^{\frac{1}{2}}} , \textbf{ }\forall f\in\mathcal{H}(\mathbb{B}), \textbf{ }\forall q\in\mathbb{B}.$
  \item[ii)]$K_{\mathcal{H}(\mathbb{B})}(q,p)=\displaystyle \sum_{k=0}^\infty Q_k(q)\overline{Q_k(p)}, \textbf{  } \forall (q,p)\in\mathbb{B}\times \mathbb{B}.$
   \item[iii)]$K_{\mathcal{H}(\mathbb{B})}(\vec{q}\,,\vec{p})=\displaystyle \sum_{k=0}^\infty (-1)^kc_k^2\vec{q}\,^k\vec{p}^k, \textbf{  } \forall (q,p)\in\mathbb{B}_0\times \mathbb{B}_0.$
   \item[iv)]$K_{\mathcal{H}(\mathbb{B})}(x,y)=\displaystyle \frac{1}{1-xy}, \textbf{  } \forall (x,y)\in ]-1,1[^2.$

\end{itemize}
\end{rem}
In the previous section we studied the notions of creation and annihilation operators associated to the Fock space in this framework. We do the same in this section for the Hardy case by studying the counterparts of the shift and backward shift operators. We keep the same definition and notation of the shift operator introduced in the expressions \eqref{Screation}, \eqref{Sexp} and Proposition \ref{shift}.
Then, we first prove the following
\begin{prop}
The shift operator $\mathcal{S}$ is a right quaternionic isometric operator from the Clifford-Appell Hardy space $\mathcal{H}(\mathbb{B})$ into itself.
\end{prop}
\begin{proof}
Let $f=\displaystyle \sum_{k=0}^\infty Q_k\alpha_k$ belongs to $\mathcal{H}(\mathbb{B})$. We apply Proposition \ref{shift} and get $$\displaystyle \mathcal{S}(f)(q)=\sum_{k=1}^\infty Q_k(q)\alpha_{k-1}, \textbf{ } \forall q\in \mathbb{B}.$$ Hence, we have

\[ \begin{split}
 \displaystyle ||\mathcal{S}(f)||^2_{\mathcal{H}(\mathbb{B})}& =\sum_{k=0}^\infty |\alpha_k|^2\\
&= ||f||^2_{\mathcal{H}(\mathbb{B})}. \\
\end{split}
\]
This shows that $\mathcal{S}$ defines an isometry on the Hardy space $\mathcal{H}(\mathbb{B})$.
\end{proof}
We will use the notation $Q_{1}^{-\odot}(q):=CK\left[\dfrac{(\vec{q})^{-1}}{c_1}\right]$ which is well defined when $\vec{q}\neq 0$. Then, in order to calculate the adjoint operator of the shift on $\mathcal{H}(\mathbb{B})$ we first prove the following result:
\begin{prop}\label{Mbshift}
For all $k\geq 1$ and $q\in\mathbb{B}\setminus \mathbb{R}$ we have $$\displaystyle(Q_{1}^{-\odot}\odot Q_k)(q)=\frac{c_k}{c_1c_{k-1}}Q_{k-1}(q), \textbf{ } $$
where $\odot$ is the C-K product and $c_l:=\displaystyle \sum_{j=0}^l(-1)^jT_j^l, \textbf{  } \forall l\geq 0.$
\end{prop}
\begin{proof}
First, we observe that $\displaystyle  (Q_{1}^{-\odot}(q))|_{x_0=0}= (Q_1(\vec{q}\,))^{-1}=\frac{(\vec{q}\,)^{-1}}{c_1}$ and $Q_{k-1}(\vec{q}\,)=c_{k-1}\vec{q}\,^{k-1}$ for $\vec{q}\neq 0$. Then, we write the series expansion associated to the C-K product and use similar techniques as we used to prove Proposition \ref{QkCKQs}. Indeed, we note that 
$$\displaystyle Q_{1}^{-\odot}(q)=\sum _{j=0}^\infty \frac{(-1)^jx_0^j}{j!}\partial_{\vec{q}\,}^j \left(\frac{(\vec{q}\,)^{-1}}{c_1} \right).$$ Moreover, for any $k\geq 1$ we have  \[ \begin{split}
 \displaystyle (Q_{1}^{-\odot}\odot Q_k)(q)& =\sum _{j=0}^\infty \frac{(-1)^jx_0^j}{j!}\partial_{\vec{q}\,}^j \left( Q_{1}^{-1}(\vec{q}\,)Q_k(\vec{q}\,) \right) \\
&= \sum _{j=0}^\infty \frac{(-1)^jx_0^j}{j!}\partial_{\vec{q}\,}^{j} \left(\frac{c_k}{c_1}\vec{q}\,^{k-1} \right)
 \\
 &= \frac{c_k}{c_1c_{k-1}}\sum _{j=0}^\infty \frac{(-1)^jx_0^j}{j!}\partial_{\vec{q}\,}^{j} \left(c_{k-1}\vec{q}\,^{k-1} \right)
 \\
 &= \frac{c_k}{c_1c_{k-1}}Q_{k-1}(q).
 \\
\end{split}
\]
\end{proof}
For all $k\geq 1$,  we introduce a family of operators defined for any $f=\displaystyle \sum_{k=1}^\infty Q_k\alpha_k $ in $\mathcal{H}(\mathbb{B})$ by \begin{equation}
\displaystyle \mathcal{M}_k(f):=\frac{c_1c_{k-1}}{c_k}Q_{1}^{-\odot}\odot f. 
\end{equation}

Then, we consider the operator obtained by applying $\mathcal{M}_k$ on each component with the corresponding degree, i.e

\begin{equation}\label{Mshift}
\mathcal{M}(f):=\displaystyle \sum_{k=1}^\infty \mathcal{M}_k(Q_k)\alpha_k.
\end{equation}
Therefore, we have an explicit expression given by
\begin{equation}\label{Mexp}
\mathcal{M}(f):=\displaystyle c_1\sum_{k=1}^\infty \frac{c_{k-1}}{c_k}[Q_{1}^{-\odot}\odot Q_k]\alpha_k.
\end{equation}
We note that using Proposition \ref{Mbshift} we can see that this operator $\mathcal{M}$ acts like the standard backward shift with respect to the Appell system $(Q_k)_{k\geq 0}$, in the sense that we have
\begin{equation}\label{Mact}
\mathcal{M}(Q_k)=Q_{k-1}, \textbf{ } \forall k\geq 1.
\end{equation}

The next result allows to compute the adjoint of the shift operator on the Hardy space $\mathcal{H}(\mathbb{B})$.
\begin{prop}
Let $f,g\in \mathcal{H}(\mathbb{B})$. Then, it holds that
$$\scal{\mathcal{M}(f),g}_{\mathcal{H}(\mathbb{B})}=\scal{f,\mathcal{S}(g)}_{\mathcal{H}(\mathbb{B})}.$$
In other words, the adjoint of the shift on $\mathcal{H}(\mathbb{B})$ is given by $$\mathcal{S}^*=\mathcal{M}.$$
\end{prop}
\begin{proof}
Let $\displaystyle f=\sum _{k=0}^\infty Q_k\alpha_k$ and $\displaystyle g=\sum _{k=0}^\infty Q_k\beta_k$ in  $\mathcal{H}(\mathbb{B})$. Thus, we have  \[ \begin{split}
 \displaystyle \mathcal{M}(f) & =\sum _{k=1}^\infty \mathcal{M}_k(Q_k)\alpha_k \\
&= \sum _{k=1}^\infty Q_{k-1} \alpha_k\\ &=\sum _{k=0}^\infty Q_{k} \alpha_{k+1}.
\end{split}
\]
We know also by Proposition \ref{shift} that
$$\mathcal{S}(g)=\sum _{k=1}^\infty Q_{k} \beta_{k-1}.$$

Therefore, we can see that  $$\scal{\displaystyle\mathcal{M}(f),g}_{\mathcal{H}(\mathbb{B})}=\sum_{k=0}^\infty \overline{\alpha_{k+1}}\beta_k=\scal{f,\mathcal{S}(g)}_{\mathcal{H}(\mathbb{B})}.$$
This ends the proof.
\end{proof}

 In \cite{AlpayShapiroVolok2013} the authors introduced a backward shift with respect to each Fueter variable using some integral operators. Inspired from this approach, we present now an equivalent way to deal with the backward shift operator in our situation. First, for all $\varepsilon>0$ we consider on $\mathcal{H}(\mathbb{B})$ a family of operators $\mathcal{R}_\varepsilon:f\longmapsto R_\varepsilon (f)$ defined using the following expression
\begin{equation}
\mathcal{R}_{\varepsilon}(f)(q):=\int_{\varepsilon}^1\frac{1}{t}\frac{\overline{\partial}}{2} \left[ f(tq)\right] dt; \textbf{ } q\in\mathbb{B}\setminus \lbrace 0 \rbrace
\end{equation}
where $\displaystyle \frac{\overline{\partial}}{2} $ denote the hypercomplex derivative with respect to the variable $q$. Then, we consider the backward shift operator given by
\begin{equation}\label{AlBas}
\mathcal{R}(f)(q):=\lim_{\varepsilon\rightarrow 0}\mathcal{R}_\varepsilon(f)(q), \textbf{ }  q\in\mathbb{B}\setminus \lbrace 0 \rbrace
\end{equation}
and \begin{equation}
\mathcal{R}f(0)=\frac{\overline{\partial}}{2}f(0).
\end{equation}

We note that the backward shift operator $\mathcal{R}$ acts by reducing the degree of the Appell system $(Q_k)_{k\geq 0}$ as follows
\begin{prop}\label{Ract}
For all $k\geq 1$, it holds that $$\mathcal{R}(Q_k)=Q_{k-1}.$$
\end{prop}
\begin{proof}
Let $k\geq 1$ and $\varepsilon>0$.
First, we note that $$Q_k(qt)=t^kQ_k(q), \textbf{  } \forall \varepsilon<t<1.$$
Then, by definition of $\mathcal{R}_\varepsilon$ and Appell property of the system $(Q_k)_{k\geq0}$ we have
  \[ \begin{split}
 \displaystyle \mathcal{R}_\varepsilon(Q_k)(q) & =\int_{\varepsilon}^1\frac{1}{t}\frac{\overline{\partial}}{2} \left[ Q_k(tq)\right] dt \\
&= \int_{\varepsilon}^1\frac{t^k}{t}\frac{\overline{\partial}}{2} \left[ Q_k(q)\right] dt\\ &=kQ_{k-1}(q)\int_{\varepsilon}^1t^{k-1}dt.
\end{split}
\]
Therefore, we obtain $$\displaystyle \mathcal{R}_\varepsilon(Q_k)(q) =Q_{k-1}(q)(1-\varepsilon^k), \forall \varepsilon>0.$$
Hence, by letting $\varepsilon\longrightarrow 0$ we conclude that $$\mathcal{R}(Q_k)=Q_{k-1}, \textbf{ }\forall k\geq 1.$$
\end{proof}
\begin{rem}
We observe thanks to formula \eqref{Mact} and Proposition \ref{Ract} that the two backward shift operators $\mathcal{M}$ and $\mathcal{R}$ coincide on the Clifford-Appell Hardy space $\mathcal{H}(\mathbb{B})$.
\end{rem}

We prove also another property related to the backward shift operator $\mathcal{R}$ on the spaces $\mathcal{HM}_b(\mathbb{B})$.

\begin{prop}\label{BSIN}
 Let  $b=(b_k)_{k\in\N}$ be a non decreasing sequence with $b_0=1$ and $f\in\mathcal{HM}_b(\mathbb{B})$. Then, the following inequality holds
\begin{equation}\label{BSin}
||\mathcal{R}(f)||_{\mathcal{HM}_b}^2\leq ||f||_{\mathcal{HM}_b}^2-|f(0)|^2.
\end{equation}

The equality holds on the Clifford-Appell Hardy space $\mathcal{H}(\mathbb{B})$.
\end{prop}
\begin{proof}
We write $f=\displaystyle \sum_{k=0}^\infty Q_k\alpha_k$ in $\mathcal{HM}_b(\mathbb{B})$. Thus, by Proposition \ref{Ract} we can see that $\displaystyle \mathcal{R}(f)=\sum_{k=0}^\infty Q_k\alpha_{k+1}.$ Therefore, using the fact that $b$ is non decreasing we get

  \[ \begin{split}
 \displaystyle ||\mathcal{R}(f)||_{\mathcal{HM}_b(\Omega)}^2& =\sum_{k=0}^\infty b_k|\alpha_{k+1}|^2\\
&\leq \sum_{k=0}^\infty b_{k+1}|\alpha_{k+1}|^2\\ &= ||f||_{\mathcal{HM}_b(\Omega)}^2-|f(0)|^2 .
\end{split}
\]
\end{proof}
\begin{rem}
We note that using Proposition \ref{BSIN} we can see that the QRKHS $\mathcal{HM}_b(\mathbb{B})$ are invariant under the backward shift $\mathcal{R}$ and they satisfy inequality \ref{BSin}. It would be intersting to investigate the relation with Schur functions and see if the converse holds also in this framework. If it is the case, it will present a counterpart of the structure result proved in Theorem 3.1.2 of \cite{AlpayRovnyakSnoo1997}.
\end{rem}
\section{The Fueter mapping range}
In this section we give an answer to Problem \ref{pb2}. Indeed, we give a characterisation of the Fueter mapping range related to the hypercomplex spaces introduced in Definition \ref{DefHb}.
\begin{thm}\label{FMR}
Let $\Omega$ be an axially symmetric slice domain and $c=(c_k)_{k\in\N}$  be a given non decreasing sequence with $c_0=1$. Then, there exists a sequence $b=(b_k)_{k\geq 0}$ such that we have $$\tau\left(\mathcal{HS}_c (\Omega)\right)=\mathcal{HM}_b(\Omega).$$
More precisely, we have
\begin{enumerate}
\item[i)] $\displaystyle b_k=\frac{c_{k+2}}{(k+1)^2(k+2)^2}, \textbf{ } \forall k\geq 0.$
\item[ii)] For all $f\in \mathcal{HS}_c(\Omega)$, we have $$||\tau(f)||_{\mathcal{HM}_b(\Omega)}=2 \sqrt{||f||_{\mathcal{HS}_c(\Omega)}^2-|f(0)|^2-c_1|f'(0)|^2}.$$
\end{enumerate}
\end{thm}
\begin{proof}
Let $g\in \tau\left(\mathcal{HS}_c (\Omega)\right), $ thus there exists $f \in\mathcal{HS}_c$ such that $g=\tau(f)$. Then, we write the series expansion  $$f(q)=\displaystyle \sum_{k=0}^\infty q^ka_k,\textbf{  } \forall q\in \Omega.$$

Thus, we have $g=\tau(f)=\displaystyle \sum_{k=0}^\infty Q_k\alpha_k, $
with $\alpha_k=-2(k+1)(k+2)a_{k+2}, \textbf{ } \forall k\geq 0.$
Now, we set $$\displaystyle b_k=\frac{c_{k+2}}{(k+1)^2(k+2)^2}, \textbf{ } \forall k\geq 0.$$
Hence, since $a_0=f(0)$ and $a_1=f'(0)$ we obtain
  \[ \begin{split}
 \displaystyle ||\tau(f)||_{\mathcal{HM}_b(\Omega)}^2& =\sum_{k=0}^\infty b_k|\alpha_k|^2\\
&= 4\sum_{k=2}^\infty c_{k}|a_{k}|^2\\ &=4\left(  ||f||_{\mathcal{HS}_c(\Omega)}^2-|f(0)|^2-c_1|f'(0)|^2 \right)<\infty.
\end{split}
\]
This ends the proof.
\end{proof}
\begin{cor}
 If we set $\mathcal{HS}_c^0:=\lbrace{f\in \mathcal{HS}_c, \textbf{ } f(0)=f'(0)=0}\rbrace$. Then, the Fueter mapping $\tau$ defines a right quaternionic isometric operator (up to constant) from $\mathcal{HS}_c^0$ onto $\mathcal{HM}_b$.
\end{cor}
\begin{proof}
We only have to apply ii) in Theorem \ref{FMR} and get
$$||\tau(f)||_{\mathcal{HM}_b(\Omega)}=2 ||f||_{\mathcal{HS}_c(\Omega)}, \textbf{  } \forall f\in \mathcal{HS}_c^0 .$$
\end{proof}
\begin{rem}
The generic calculations provided in Theorem \ref{FMR} confirm the results obtained in \cite{DKS2019} for the Fock and Bergman cases.

\end{rem}
\begin{rem}
We note that in Theorem \ref{FMR} even if the sequence $b$ is not necessarily a non decreasing sequence but the corresponding spaces $\mathcal{HM}_b$ are QRKHS. For the Fock-Fueter space on $\Hq$ we refer to the calculation details provided in \cite{DKS2019}. However, on the quaternionic unit ball $\mathbb{B}$ this fact results thanks to the convergence of a certain power series associated to the sequence $b$.
\end{rem}

\begin{prop}  Let $c$ and $b$ two sequences as in Theorem \ref{FMR}. Then, the power series given by

\begin{equation}
\displaystyle \sum_{k=0}^\infty \frac{|q|^{2k}}{b_k}=\sum_{k=}^{\infty}\frac{(k+1)^2(k+2)^2}{c_{k+2}}|q|^{2k},
\end{equation}

is convergent on the quaternionic unit ball $\mathbb{B}$.
\end{prop}
\begin{proof}
Let $q\in\mathbb{B}$ and set $$s_k=\frac{(k+1)^2(k+2)^2}{c_{k+2}}|q|^{2k}, \textbf{  } \forall k\geq 0.$$
We have $$\displaystyle \frac{s_{k+1}}{s_k}=|q|^2\frac{(k+3)^2c_{k+2}}{(k+1)^2c_{k+3}}, \textbf{ } \forall k\geq 0.$$

Then, using the fact that the sequence $(c_k)_{k\geq 0}$ is non decreasing we can see that $$\lim_{k\rightarrow \infty}\frac{s_{k+1}}{s_{k}}\leq |q|^2<1.$$
Hence, by the d'Alembert ratio test the thesis follows.
\end{proof}
\begin{rem}
As a consequence of the previous Proposition it is not difficult to see that on $\mathbb{B}$ the hypercomplex space $\mathcal{HM}_b$ obtained in Theorem \ref{FMR} is a QRKHS with a reproducing kernel given by
\begin{equation}
K_{\mathcal{HM}_b}(q,p)=\displaystyle \sum_{k=0}^{\infty}\frac{(k+1)^2(k+2)^2}{c_{k+2}}Q_k(q)\overline{Q_k(p)}, \textbf{  } \forall (q,p)\in \mathbb{B}\times \mathbb{B}.
\end{equation}
\end{rem}
In the following table we list some spaces of slice hyperholomorphic functions and their Fueter mapping ranges denoted respectively by $\mathcal{HS}_c$ and $\mathcal{HM}_b$, the associated sequences $c$ and $b$ and the Fueter mapping norms.

\begin{table}[ht]
\caption{Some spaces $\mathcal{HM}_b$ obtained in Theorem \ref{FMR}} 
\centering 
\begin{tabular}{c c c c} 
\hline\hline 
 $\mathcal{HS}_c$  & $c_k$ & $b_k$ & $||\tau(f)||_{\mathcal{HM}_b}$  \\ [0.5ex] 
\hline 
Hardy & $1$ & $\dfrac{1}{(k+1)^2(k+2)^2}$ & $2 \sqrt{||f||_{\mathcal{HS}_c}^2-|f(0)|^2-|f'(0)|^2}$ \\ 
Fock & $k!$ & $\dfrac{k!}{(k+1)(k+2)}$ & $2 \sqrt{||f||_{\mathcal{HS}_c}^2-|f(0)|^2-|f'(0)|^2}$ \\
Dirichlet & $k$ & $\dfrac{1}{(k+1)^2(k+2)}$ & $2 \sqrt{||f||_{\mathcal{HS}_c}^2-|f(0)|^2-|f'(0)|^2}$ \\
Bergman & $\dfrac{1}{k+1}$ & $\dfrac{1}{(k+3)(k+1)^2(k+2)^2}$ & $2 \sqrt{||f||_{\mathcal{HS}_c}^2-|f(0)|^2-\frac{1}{2}|f'(0)|^2}$ \\
 [1ex] 
\hline 
\end{tabular}
\label{table:nonlin} 
\end{table}

\nocite{*}

\end{document}